\documentclass[12pt]{amsart}
\usepackage{amssymb}
\usepackage{amsfonts}
\usepackage{amsthm}
\numberwithin{equation}{section}
\theoremstyle{plain}
\newtheorem{theorem}{Theorem}[section]
\newtheorem{observation}[theorem]{Observation}
\newtheorem{proposition}[theorem]{Proposition}
\newtheorem{corollary}[theorem]{Corollary}

\theoremstyle{definition}
\newtheorem{definition}[theorem]{Definition}
\theoremstyle{remark}
\newtheorem*{remark}{Remark}
\newtheorem*{example}{Example}
\newcommand{\refE}[1]{(\ref{E:#1})}
\newcommand{\refS}[1]{Section~\ref{S:#1}}
\newcommand{\refSS}[1]{Section~\ref{SS:#1}}
\newcommand{\refT}[1]{Theorem~\ref{T:#1}}
\newcommand{\refO}[1]{Observation~\ref{O:#1}}
\newcommand{\refP}[1]{Proposition~\ref{P:#1}}
\newcommand{\refD}[1]{Definition~\ref{D:#1}}
\newcommand{\refC}[1]{Corollary~\ref{C:#1}}

\newcommand{\R}{\ensuremath{\mathbb{R}}}
\newcommand{\C}{\ensuremath{\mathbb{C}}}
\renewcommand{\S}{\ensuremath{\Sigma}}
\newcommand{\N}{\ensuremath{\mathbb{N}}}

\renewcommand{\P}{\ensuremath{\mathbb{P}}}
\newcommand{\Z}{\ensuremath{\mathbb{Z}}}

\newcommand{\Hi}{\ensuremath{\mathcal{H}}}

\newcommand{\nl}{\hfill\newline}
\renewcommand{\i}{{\,\mathrm{i}\,}}

\newcommand{\e}{\mathrm{e}}
\newcommand{\Ls}{\mathrm{L}^2(M,L)}
\newcommand{\Hop}[1][m]{\mathrm{H}^0(M,L^{#1})}
\newcommand{\ghm}[1][m]{\Gamma_{hol}(M,L^{#1})}
\newcommand{\gh}{\Gamma_{hol}(M,L)}
\newcommand{\ghb}{\Gamma_{hol}^b(M,L)}

\newcommand{\gul}{\Gamma_{\infty}(M,L)}

\renewcommand{\d}{\partial}
\newcommand{\db}{\overline{\partial}}
\newcommand{\wb}{\overline{w}}
\newcommand{\zb}{\overline{z}}

\newcommand{\qlb}{\ensuremath{(L,h,\nabla)}}
\newcommand{\w}{\ensuremath{\omega}}
\newcommand{\pnc}[1][N]{\ensuremath{\P^{#1}(\C)}}
\newcommand{\PGL}{\mathrm{PGL}}
\newcommand{\PU}{\mathrm{PU}}
\newcommand{\volu}{\mathrm{vol}}
\newcommand{\skp}[2]{{\langle #1,#2\rangle}}
\newcommand{\skps}[3]{{\langle #1,#2\rangle}_{#3}}
\newcommand{\divi}[1]{{\langle #1\rangle}}
\begin{document}
\vspace*{-1cm}
{\hspace*{\fill} Mannheimer Manuskripte 232}

{\hspace*{\fill} math.DG/9903105}

\vspace*{2cm}

\title{Coherent state embeddings, polar divisors and Cauchy formulas}
\thanks{This work was partially supported by the Volkswagen-Stiftung
(RiP-program at Oberwolfach)}
\author{Stefan Berceanu}
\address[Stefan Berceanu]{Institute for Physics and Nuclear Engineering\\
         Department of Theoretical Physics\\
         PO BOX MG-6, Bucharest-Magurele\\
         Romania}
\email{Berceanu@theor1.theory.nipne.ro}
\author{Martin Schlichenmaier}
\address[Martin Schlichenmaier]{Department of Mathematics and 
  Computer Science, University of Mannheim, D7, 27 \\
         D-68131 Mannheim \\
         Germany}
\email{schlichenmaier@math.uni-mannheim.de}
\begin{abstract}
For arbitrary quantizable compact K\"ahler manifolds,
relations between
the geometry given by the coherent states based on the 
manifold  and the 
algebraic (projective) geometry realised   via the coherent state mapping
into 
projective space,    are studied.
Polar divisors, formulas relating the scalar products
of coherent vectors on the manifold with the corresponding 
scalar products on  projective  space
(Cauchy formulas), two-point, three-point and more 
generally cyclic $m$-point functions are discussed.
The three-point function is related to 
the shape invariant of geodesic triangles
in projective space.
\end{abstract}
\date{17.3.99}
\maketitle
\section{Introduction}\label{S:intro}
In this article the close
relations between the coherent state approach appearing in
quantum mechanics and certain aspects of algebraic geometry,
respectively K\"ahler geometry are considered.
We analyse the case where the phase-space manifold of
the theory is a compact K\"ahler manifold $(M,\w)$.
The symplectic structure which gives the kinematics
of the theory is defined via the K\"ahler form $\w$.
The geometric quantization condition requires the existence of a 
line bundle (i.e. the quantum line bundle) with curvature essentially
equal to the K\"ahler form.
This implies that the phase-space manifold is projective
algebraic. Hence it admits an embedding 
into projective space.
An embedding can be explicitly given by the global sections
of a suitable tensor power of the quantum line bundle.
It is usually known as Kodaira embedding.
Vice versa, every submanifold of projective  
space is a quantizable K\"ahler manifold.

Berezin's coherent states \cite{Berebd}
in their reformulation and generalization due to Rawnsley \cite{Raw}
define also an embedding into projective space.
It turns out that this embedding is nothing else 
as the Kodaira embedding (respectively its 
conjugate) with respect to an orthonormal basis
of the space of global holomorphic sections 
of a suitable tensor power of the quantum line bundle.
Here  the scalar product is induced by the K\"ahler form
(see \refS{coher} for details).
Having the possible inhomogeneous readership in mind we recall the 
basics of the afore mentioned concepts in \refS{coher}.
In this way we also explain our notation and prove some results
used later on.

The main goal of this article is to study relations between the geometry of
the quantizable K\"ahler manifold using  coherent states and the 
algebraic (or projective) geometry of the embedded manifold 
in $\pnc$.
Such kind of relations for homogeneous manifolds 
(with respect to Perelomov's coherent states) 
were studied by Berceanu
 \cite{Berpre94},\cite{BerJGP},\cite{Berbia96},\cite{Berbia98},\cite{Berbia97}.
Here we will make similar definitions and 
prove some analogous results
for arbitrary compact K\"ahler manifolds. Clearly,  now we have to use 
the coherent states of Berezin-Rawnsley.

The first objects we introduce are the polar divisors (\refS{polar}).
 The polar divisor of
a point on $M$ is the divisor consisting of the points on the manifold
whose coherent states are orthogonal to the  coherent state
associated to the fixed point.
It turns out that the polar divisor is indeed a divisor in the sense of
algebraic geometry.
It should not be confused with the divisor of the
polar part of a meromorphic function.

The polar divisors are useful for many purposes.
It was shown by Berceanu 
that for  Grassmannians and more general for symmetric
spaces 
\cite{BerJGP},\cite{BerGr},\cite{Berbia96}
the polar divisor $\Sigma_x$ (with respect to Perelomov's
coherent states)  coincides with the cut-locus of the point $x$.
For general compact K\"ahler manifolds
the polar divisors  describe
the zero-sets 
of two-point functions (and via them also of the $m$-point functions).
They  appear as singularity sets of the analytic extensions
for real-analytic metrics in the bundle 
and as  singularity sets of the covariant two-point
Berezin symbols.

Next (\refS{cauchy})
we discuss ``Cauchy formulas''. Under a "Cauchy formula" we understand
a relation between the scalar product of the
coherent states (more precisely, of the
coherent vectors) associated to two points on the manifold 
(again more precisely, associated to points in the total space
of the quantum line bundle) 
and the 
scalar product of the via the coherent state map embedded two points
(more precisely, the scalar product of certain homogeneous
representatives of the embedded points).
The main results  are contained in  \refT{cauchy} and the propositions
in  \refS{cauchy}.
The denomination ``Cauchy formula'' was used 
 in this context the first time in \cite{Berpre94}
for Perelomov's coherent states on flag manifolds. For the Grassmannian
 the  appearing formulas are 
essentially the (Binet-)Cauchy formulas \cite[p.10]{Gant}
which give relations between the intrinsic metric on the
Grassmannian and the pull-back of the Fubini-Study metric 
obtained via the Pl\"ucker embedding.

In the remaining part of \refS{cauchy} the 
two-point function  and cyclic $m$-point functions are discussed.
Considered on $M\times M$
the complex-valued two-point function 
has a phase ambiguity.
This ambiguity 
can  either be removed 
by considering the modulus of the function 
or alternatively by fixing a holomorphic 
section of the quantum line bundle as a reference lift
to the quantum line bundle.
In the first case one ends up with the two-point function 
studied by Cahen, Gutt and Rawnsley \cite{CGR2} which is
related to Calabi's diastatic function.
But it turns out that the complex-valued ``non-canonical'' two-point function
plays at least a very useful intermediate r\^ole.
The polar divisors appear naturally 
in this context as zero sets of  two-point functions.
Next, cyclic $m$-point functions are introduced.
They are well-defined on $M^{\times m}$ and invariant under cyclic
permutations of its arguments.
It is shown that they are invariant under pull-back via the coherent state 
embedding into projective space.
They can be expressed in terms of the 
 Cayley distances  of the 
 embedded  points in 
 projective space  and a phase factor, depending on the points.
The three-point function is studied in more detail.
Here the phase is related 
to  the shape
invariant of the geodesic triangle which has the embedded points as vertices.
The shape invariant was introduced 1939 by Blaschke and Terheggen \cite{BlTe}.
By a result of Hangan and Masala \cite{HaMa} the phase can be calculated
via integrating the K\"ahler form over geodesic triangles.
See the closing \refT{geo} for the detailed result.    
\section{Coherent state embedding}\label{S:coher}
\bigskip

\subsection{Quantizable K\"ahler manifolds and Kodaira embedding\hfill}
\label{SS:kodaira}

{~}  
\vspace*{0.2cm}

Let $(M,\w)$ be a K\"ahler manifold of complex dimension $n$,
 i.e. $M$ a complex manifold and
$\w$ a K\"ahler form on $M$. In the following we will mainly consider
compact K\"ahler manifolds. If nothing else is said 
we will assume compactness.
A further data we need is the triple $\qlb$, with  a holomorphic
line bundle  $L$ on $M$, a hermitian metric $h$ on $L$ 
(with the convention that
it is conjugate linear in the first argument) and a connection
$\nabla$ compatible with the metric on $L$ and the complex structure.
 With respect to local holomorphic coordinates of the manifold
and with respect to a local holomorphic frame for the bundle
the metric $h$ can be given as 
\begin{equation}\label{E:locmet}
 h(s_1,s_2)(x)=\hat h(x)\overline{\hat s_1}(x)\hat s_2(x),
\end{equation}
where $\hat s_i$ is a local representing function for the
section $s_i$ ($i=1,2$) and $\hat h$ is a locally defined real-valued 
function on $M$.
The compatible connection is  uniquely defined and is given
in the local coordinates as
$\ \nabla=\d +(\d\log \hat h) +\db$.
The curvature of $L$ is defined as the two-form 
\begin{equation}\label{E:curvdef}
F(X,Y)=\nabla_X\nabla_Y-\nabla_Y\nabla_X-\nabla_{[X,Y]}\ ,
\end{equation}
where $X$ and $Y$ are vector fields on $M$.
In the local coordinates the curvature can be expressed as
\begin{equation}\label{E:curvloc}
F=\db\d\log \hat h=-\d\db\log \hat h\ .
\end{equation}

A K\"ahler manifold $(M,\w)$ is called \emph{quantizable}
if there
exists  such a triple $\qlb$ which obeys
\begin{equation}\label{E:quantcond}
F(X,Y)=-\i\w(X,Y)\ .
\end{equation}
The condition \refE{quantcond} is called the (pre)quantum condition.
The bundle $\ \qlb\ $ is called a (pre)quantum line  bundle.
Usually we will drop $\nabla$ and sometimes also
$h$  in the notation.

For the following we assume  $(M,\w)$  to be a quantizable
K\"ahler manifold with quantum line bundle  $\qlb$.
Let us note some important consequences of the quantum condition.
Firstly, 
 we get for  the Chern form of the line bundle $L$
the relation
\begin{equation}\label{E:chernf}
c(L)=\frac {\i}{2\pi}F=\frac {\w}{2\pi}\ .
\end{equation}
This implies that $L$ is a positive line bundle. In the terminology
of algebraic geometry it is an ample line bundle.
This says that there exists a tensor power $L^{m_0}:=L^{\otimes m_0}$
with $m_0$ a positive integer such that $M$ can be holomorphically
embedded into projective space $\pnc[N]$ 
using the global holomorphic sections of
$L^{m_0}$.
Let us describe this embedding in more detail.
We will denote the space of global holomorphic sections
by $\Hop[{m_0}]$, or depending on the context, also as $\ghm[{m_0}]$
(if we regard it as subspace of the space of differentiable sections),
resp. by $\Hi$ (if we regard it as the quantum Hilbert space).
By  compactness of the manifold $M$ this vector space is 
finite-dimensional. We take 
 $N=\dim\Hop[{m_0}]-1$ and after fixing a basis of the global 
sections the embedding is given as
\begin{equation}\label{E:embedd}
\varphi:M\ \hookrightarrow\  \pnc[N],
\qquad z\ \mapsto\  \varphi(z)=(s_0(z):s_1(z):\ldots:s_N(z))\in \P^N(\C)\ .
\end{equation}
Here we denote the point $\varphi(z)$ in projective space by its
homogeneous coordinates.
Recall that  two sets of homogeneous coordinates correspond to the same point
if
and only if they are a non-zero scalar multiple of each other.
To evaluate the sections one chooses local representing functions
for the sections. 
Clearly, they are only well-defined up to a common scalar function.
Hence, only after  passing to the projective space the map will be 
well-defined.
The conclusion that from the positivity of the line bundle 
it follows that there exists such an embedding is the 
content of 
Kodaira's embedding theorem, see \cite{GH}, \cite{Wells}.
By Chow's theorem \cite[p.166]{GH} compact 
submanifolds of $\pnc[N]$ are projective
varieties, i.e. they can be given as zero-sets of a finite number
of homogeneous polynomials in the coordinates of  $\pnc[N]$,
see \cite{SchlRS}, \cite{GH}.
Hence, we obtain the first part of the following important
observation
\begin{observation}\label{O:proj}
Quantizable compact K\"ahler manifolds are submanifolds of 
$\pnc[N]$, hence projective algebraic. Conversely, every projective
algebraic manifold will be a quantizable K\"ahler manifold.
\end{observation}
The second part will follow from the discussion further down in this
section.

In the language of K\"ahler geometry
quantizable compact K\"ahler manifolds are Hodge manifolds \cite{Wells}, 
\cite{GH}.
This is due to the fact that from the relation \refE{chernf} 
it follows that the class of $\w$ is $2\pi$ times  a Chern class, hence 
a integral class, i.e. a class which gives an integer
when integrated over a closed 2-surface in $M$.

 The number $N$ can
be explicitly given with the help of the
Grothendieck-Hirzebruch-Riemann-Roch
Theorem \cite{Hart}, \cite{SchlRS}.
A different choice of basis taken for the
embedding corresponds to a holomorphic automorphism
of
$\pnc[N]$, i.e. to an element of 
\nl $\PGL(N+1,\C)$
mapping the images onto each other.

In the following we will assume that $L$ is already
very ample.
This says that $M$ can already be embedded using the global sections of $L$.
 If this is not yet the case we  
can always choose a $m_0\in\N$ such that $L^{m_0}$ is very ample.
Now in  generality  for a $m^{th}$ tensor power 
$L^m$ of the bundle $L$ a metric $h^{(m)}$ 
and a connection $\nabla^{(m)}$ is given by
\begin{equation}
\begin{gathered}
h^{(m)}:=\underbrace{h\otimes\cdots\otimes h}_{\text{$m$ times}},\\
\nabla^{(m)}=
\nabla\otimes 1\otimes\cdots\otimes 1+
1\otimes\nabla\otimes \cdots\otimes 1+\cdots+
1\otimes\cdots\otimes 1\otimes\nabla\ .
\end{gathered}
\end{equation}
The corresponding local objects are 
(if one takes the $m^{th}$ tensor power of the frame of $L$ as frame 
for  $L^{m}$)
\begin{equation}
\widehat{h^{(m)}}=(\hat h)^m,\qquad
\nabla^{(m)}=\partial+m\,(\partial\log \hat h) +\db,\qquad
F^{(m)}=mF=-\i m\,\w\ .
\end{equation}
Hence, for every $m\in\N$ the bundle  $\ (L^m,h^{(m)},\nabla^{(m)})\ $ 
is a quantum line bundle for the K\"ahler manifold
 $(M,m\,\w)$.
Note that  the underlying complex manifold remains the same, only
the K\"ahler form is multiplied
by an integer\footnote
{The process that starting from one line bundle $L$ one obtains for 
every $m\in\N$ a quantization allows to introduce semi-classical limits of
the quantization scheme (geometric quantization, Berezin-Toeplitz
quantization,
coherent state quantization,..), to prove approximation results
for them (e.g. see \cite{BMS}), and to show the existence of star products
\cite{CGR1},\cite{CGR2},\cite{Schlbia95},\cite{Schldef}.}.
So, if we start with $(M,m_0\,\omega)$ the corresponding quantum line bundle 
$L^{m_0}$ is very ample.

\medskip
A second consequence of the quantum condition \refE{quantcond}
is that the metric in the quantum bundle can be expressed with the
help of a local K\"ahler potential \cite{ChCM}. For a K\"ahler
manifold
there  exist
locally real-valued (non-unique) functions $K$ such that 
$\w=\i\d\db K$.
With  the quantum condition \refE{quantcond} it follows from \refE{curvloc}
\begin{equation}\label{E:kaehlerloc}
\w=\i\db\d\log\hat h\ .
\end{equation}
Hence a local K\"ahler potential can be given as 
\begin{equation}\label{E:kaehlerpot}
K(z)=-\log\hat h(z),\quad\text{resp.},\quad
\hat h(z)=\exp(-K(z))\ .
\end{equation}

\bigskip
Recall the K\"ahler structure of the projective space $\pnc$.
The points $[z]$ in  $\pnc$ are given by their homogeneous coordinates
$\ [z]:=(z_0:z_1:\ldots:z_N)\ $. 
In the affine chart $V_0$ consisting
of the points with $z_0\ne 0$ we
take $\ w_j=z_j/z_0\ $ with $j=1,\ldots, N$ as holomorphic
coordinates.
In the similar way we define affine charts $V_k$, $k=1,\ldots,N$
and corresponding holomorphic coordinates.
The union $\bigcup_{k=0}^N V_k$ is now an affine covering of $\pnc$.
Denote by 
\begin{equation}\label{E:proj}
\tau:\ \C^{N+1}-\{0\}\quad \to\quad \pnc
\end{equation} 
the projection which is obtained by identifying the
whole line through
$0$ and the point z with the point in 
projective space with  homogeneous coordinates $[z]$.
The K\"ahler form on \pnc\ is the Fubini-Study fundamental
form. On $V_0$ it is given as
\begin{equation}\label{E:fubstud}
\w_{FS}:=
\i\frac
{(1+{\Vert w\Vert}^2)\sum_{i=1}^N dw_i\wedge
d\wb_i-\sum_{i,j=1}^N\wb_iw_jdw_i\wedge d\wb_j} {{(1+{\Vert w\Vert}^2)}^2}
\ .
\end{equation}
Here ${\Vert w\Vert}^2:=\sum_{i=1}^N \wb_i w_i$, as usual.
Alternatively it can be described as \cite{Kob}
\begin{equation}\label{E:fubstudhom}
\tau^*\w_{FS}(z) =\i\db\d\log {\Vert z\Vert}^2\ .
\end{equation} 
Over $\pnc$
we have  the tautological line bundle $U$.
Its fiber over $[z]$ consists of the line 
through $0$ and $z$.
Taking the standard metric in $\C^{N+1}$ it is endowed
with a natural hermitian fiber 
metric. Note that the manifold
$\C^{N+1}-\{0\}$ can be identified with the total space of $U$ with
the zero section removed.
With respect to the affine chart $V_0\cong\C^N$ we can write two
elements of the same fiber over $w\in\C^N$ as
\begin{equation} 
s_1=\alpha\cdot(1,w_1,\ldots,w_N),\qquad
s_2=\beta\cdot(1,w_1,\ldots,w_N),\qquad
\end{equation}
and  obtain
\begin{equation}
\bar s_1\cdot s_2=\bar\alpha\cdot \beta\cdot (1+{\Vert w \Vert}^{2})\ .
\end{equation}
Hence, the local representing function  in $V_0$  for the hermitian metric
of the line bundle $U$  
(and with respect to the standard frame
$V_0\to \C\times V_0,\ w\mapsto (1,w)$)
is 
\begin{equation}\label{E:tautomet}
\hat k(w)=1+{\Vert w \Vert}^{2}\ .
\end{equation}

The quantum line bundle is the dual of the 
tautological bundle, the hyperplane bundle $H=U^*$.
The hermitian metric of the hyperplane bundle
can be given in the affine chart by the representing 
function
\begin{equation}\label{E:hypmet}
\hat h(w)=\frac {1}{1+{\Vert w \Vert}^{2}}\ .
\end{equation}
The  global holomorphic sections of $H$ can be identified with
the linear forms in the $N+1$ coordinate functions  $Z_i$.

We were using the term "the quantum line bundle" indicating
that there is up to algebraic isomorphy just one line bundle with 
curvature form  $-\i\w$. In general this is not the case.
But for the projective space 
there is for every degree up to isomorphy just one line bundle 
and the degree is fixed 
by the curvature.
Hence for the projective space (with the 
Fubini-Study K\"ahler form) the quantum line bundle is fixed. 
In fact the same is true for any simply-connected, compact, 
quantizable K\"ahler manifold.
In this case there is at most one line bundle which has a given candidate 
as curvature form, see \cite[Thm. 2.2.1]{Kost}.
Here a warning is in order.
It is not excluded that for the same underlying complex manifold
there exist (essentially) different K\"ahler forms
and hence essentially different associated quantum line bundles.

\medskip
If $M$ is a projective submanifold of $\pnc$ with 
\begin{equation}\label{E:embed}
i:M\hookrightarrow \pnc
\end{equation}
the inclusion then $(M,i^*\w_{FS})$ is a K\"ahler manifold
\cite{GH}, \cite{Wells} which 
is quantizable with the associated quantum line bundle
$(i^*H,i^*h)$. Here $i^*$ is nothing else as the restriction to 
the submanifold.
Hence, projective manifolds are quantizable K\"ahler manifolds.
This shows the second statement in \refO{proj}.

We have to stress an important fact.
If $(M,\w_M)$ is a quantizable K\"ahler manifold with very ample quantum line
bundle $L$ then we saw that $L$ induces an embedding
$i:M\hookrightarrow \pnc$. By the construction $i^*H\cong L$ as holomorphic
line bundle. Now $(M,i^*\w_{FS})$ is a K\"ahler manifold with 
the same underlying complex manifold structure. But in general
$\w_{M}\ne i^*\w_{FS}$, so the  K\"ahler structure
of $M$ does  not coincide  with the induced K\"ahler structure 
coming from
the embedding.
The embedding is in general not an isometric (K\"ahler) embedding.
The situation is very much related to Calabi's diastatic function,
\cite{Cal}, \cite{CGR1}.
In general we only know the identity of the deRham classes
$\ [\w_M]=[i^*\w_{FS}]$.
This follows from the identity of the Chern classes
$c_1(L)=c_1(i^*H)=i^*c_1(H)$ and from the fact that
by the quantum condition the K\"ahler forms represent
(up to a factor) the curvature class.

In the compact case
the K\"ahler form $\w$ 
(hence the metric on $M$) fixes via the quantum condition the hermitian 
bundle metric in $L$ up to a scalar constant:
\begin{proposition}\label{P:invarmet}
Let $(M,\w)$ be a quantizable compact K\"ahler manifold with
quantum bundles $(L,h)$ and $(L,h')$ then
$h=\mathrm{e}^{\alpha} \cdot h'$, with 
$\alpha\in \R$.
\end{proposition}
\begin{proof}
Represent the metric $h$ and $h'$ with respect to a local frame of
the (same) bundle $L$ as local functions $\hat h$ and $\hat h'$.
By the quantum condition \refE{quantcond}
we obtain for the K\"ahler form
\begin{equation}
\w=\i\d\db\log \hat h=
\i\d\db\log \hat h'.
\end{equation}
Hence $\d\db(\log \hat h-\log \hat h')=0$ or equivalently
$\log ({\hat h}/{\hat h'})$ is a locally defined harmonic function.
But the quotient of the two metrics is a globally defined function.
Hence $\log ({\hat h}/{\hat h'})$ 
 is a globally defined harmonic function on the compact
manifold and hence a constant $\alpha\in\R$.
This shows the claim.
\end{proof}

\subsection{Embedding via coherent states}\label{SS:coherent}
\hphantom{hallo}
{~}  

\vspace*{0.2cm}
We now want  to describe an (anti-)holomorphic embedding of the
K\"ahler manifold $(M,\w)$ into projective space 
using  coherent states.
We use Berezin's coherent states \cite{Berebd}, \cite{Berequ}, \cite{Beress}
in the coordinate independent global
version due to Rawnsley, see \cite{Raw}, \cite{CGR1}.

First we have to introduce a scalar product in the space of global
holomorphic sections of the quantum line bundle $L$.
With the normalized volume form
\begin{equation}
\Omega:=(-1)^{\binom {n}{2}}\frac {1}{n!}\;
\underbrace{\w\wedge\ldots\wedge\w}_{\text{$n$ times}}\ 
\end{equation}
and with the fiber metric $h$ we can introduce a scalar product and a
norm on the space of differentiable sections $\gul$
\begin{equation}\label{E:skp}
\langle\varphi,\psi\rangle:=\int_M h(\varphi,\psi)\;\Omega\  ,
\qquad
||\varphi||:=\sqrt{\langle \varphi,\varphi\rangle}\ .
\end{equation}
Using local representing functions $\hat\varphi$ and $\hat\psi$
for the sections and $\hat h$ for the metric the scalar product can
be described as 
\begin{equation}
\langle\varphi,\psi\rangle=\int_M \hat h(z)\,\overline{\hat\varphi(z)}
\,\hat\psi(z)\,\Omega(z)=
\int_M \exp(-K(z))\,\overline{\hat\varphi(z)}\,\hat\psi(z)\,\Omega(z)
\ .
\end{equation}
In the second form we used the local K\"ahler potential \refE{kaehlerpot}.
Clearly, these integrals should be calculated  locally and 
their values patched together by a partition of unity argument\footnote
{
Sometimes it is useful to write $\hat h(z,\zb)$, resp.
$K(z,\zb)$ to remind of the fact that these functions are
not holomorphic in $z$ and      
(even more important) to consider the possibility to extend the objects
analytically to $h(z,\bar w)$, resp. $K(z,\bar w)$. }.

The scalar product can be restricted to the finite dimensional
subspace of global holomorphic sections.

\medskip
Recall that we assume  the quantum line bundle to be  already very ample.
Denote by $\pi:L\to M$ the bundle projection and by $L_0$ the total space of
$L$ with the zero section $0(M)$ removed.
Fix $q\in L_0$ and take an arbitrary holomorphic
section $s$ of $L$. By evaluation of the section 
at $x=\pi(q)$ the relation
\begin{equation}\label{E:qhdef}
s(\pi(q))=\hat q(s)\cdot q \
\end{equation}
defines a linear form 
\begin{equation}
\hat q\ :\ \gh\ \to\ \C,\qquad
s\mapsto\hat q(s)\ .
\end{equation}
Using the scalar product on the space of global sections, 
by Riesz's theorem there exists exactly one holomorphic section 
$e_q$ with 
\begin{equation}\label{E:qhdefa}
\langle e_q,s\rangle=\hat q(s),\qquad\text{for all }\quad s\in\gh\ .
\end{equation}
If we choose an orthonormal basis 
 $s_j,\ j=0,\ldots,N:=\dim\gh-1$ then $e_q$ can be explicitly 
given as 
\begin{equation}\label{E:cohdecomp}
e_q=\sum_{j=0}^N\overline{\hat q(s_j)}\,s_j\ .
\end{equation}
Let $x=\pi(q)$ and choose $q'\in \pi^{-1}(x)$ with $q'\ne 0$ then
there is a $c\in\C^*$ with $q'=cq$. 
{}From \refE{qhdef} we conclude
$\ \widehat{q'}=c^{-1}\widehat{q}$ and using \refE{qhdefa}
we obtain
\begin{equation}\label{E:cohtrans}
e_{cq}=\bar c^{-1}\cdot e_q\ .
\end{equation}
We obtain two mappings
\begin{alignat} {3}
&L_0\to\gh^*,&\qquad& q\mapsto \hat q\ ,&\qquad&\text{and}
\\ \label{E:cohema}
&L_0\to\gh,&\qquad& q\mapsto e_q\ .
\end{alignat}
The first one is holomorphic, the second one antiholomorphic.
By the above relations both maps are well-defined on $M$ if we pass to
the
projectivized vector spaces
\begin{alignat} {3}
&M\to\P(\gh^*), &\qquad& x\mapsto [\widehat{\pi^{-1}(x)}]\ ,&\qquad&\text{and}
\\ \label{E:cohemap}
&M\to\P(\gh), &\qquad&x\mapsto [e_{\pi^{-1}(x)}]\ .
\end{alignat}
Here $[v]$ denotes the equivalence class of a vector $v$ of a vector space
$V$ in the
projectivized  vector space $\P(V)$.
In abuse of notation we understand by $\pi^{-1}(x)$ only the 
non-zero elements of the fiber over $x$.

Note that  $\hat q\equiv 0$ or equivalently $e_q\equiv 0$ would imply  
that all sections $s\in\gh$ will vanish at $\pi(q)$ 
and this  contradicts the very ampleness of $L$.

Depending on $q\in L_0$ the sections $e_q\in\gh$ are called
\emph{coherent vectors}.
Depending on $x\in M$ the $ [e_{\pi^{-1}(x)}]\in\P(\gh)$
are called \emph{coherent states}.
To simplify the notation we will set $\ e_x:= [e_{\pi^{-1}(x)}]$.
The  mappings \refE{cohema} and \refE{cohemap} 
are the \emph{coherent vector mapping}, resp.
the \emph{coherent state mapping}.

\medskip
To identify $\P(\gh)$ with $\pnc$ we  choose an 
orthonormal basis. The description \refE{cohdecomp} shows that the
coherent state mapping is given as 
\begin{equation}\label{E:cohmap}
x\ \mapsto\  e_x=[e_{\pi^{-1}}(x)]\ \mapsto\quad
(\overline{\hat q(s_0)}:\overline{\hat q(s_1)}:...:
\overline{\hat q(s_N)})
=
(\overline{s_0(x)}:\overline{s_1(x)}:...:
\overline{s_N(x)})\ ,
\end{equation}
For the last equality we used
$s_j(x)=s_j(\pi(q))=\hat q(s_j)\cdot q$.
\begin{proposition}\label{P:cohemb}
The map \refE{cohmap}
\begin{equation}\label{E:cohemb}
M\ \to\ \P(\gh)\cong\pnc[N],
\end{equation}
is an antiholomorphic embedding.
Up to complex conjugation it coincides with the Kodaira embedding
\refE{embedd} obtained with respect to the chosen orthonormal basis.
\end{proposition}
\begin{proof}
That the map is well-defined we showed above. That it is an embedding
follows
from the  observation  that
Equation \refE{cohmap} is up to 
complex conjugation nothing else as the Kodaira embedding with
respect to the very ample line bundle $L$.
\end{proof}
In the following it will be more convenient to consider the 
complex conjugate of the coherent state embedding
\refE{cohmap}
\begin{equation}\label{E:cohmapc}
x\ \mapsto\  \overline{e_x}=[\overline{e_{\pi^{-1}}(x)}]\ \mapsto\quad
({\hat q(s_0)}:{\hat q(s_1)}:...:
{\hat q(s_N)})
=
({s_0(x)}:{s_1(x)}:...:
{s_N(x)})\ ,
\end{equation}
which is a holomorphic embedding. We will use 
the term {\it coherent state embedding} also for \refE{cohmapc}
if there is no danger of confusion.

Note that a different orthonormal basis (ONB)
 will yield an embedding which is equivalent
under a $\PU(N+1)$ action to the chosen one.

In the language of physics  \refP{cohemb} means that the
phase space of a mechanical system (assumed here to be
K\"ahlerian) can be embedded via coherent states into a 
projectivized Hilbert space, the quantum Hilbert space.

It should be pointed out that the coherent state embedding is not just
Kodaira embedding. It is Kodaira embedding using orthonormal sections.
The scalar product used to define the orthonormality on $\gh$ (which 
should be interpreted as the quantum Hilbert space $\mathcal{H}$)
 is induced by 
the K\"ahler form on the manifold and by the hermitian metric in the
bundle.
In view of the quantization condition
the latter itself can be 
related to the K\"ahler form of the manifold,
see \refE{kaehlerloc} and \refP{invarmet}.
The K\"ahler form (interpreted as symplectic form) is an important
ingredients to the description of the system to be quantized.

\medskip
If one considers non-compact K\"ahler manifolds then the scalar
product \refE{skp} on $\gul$ 
or more precisely on $\Ls$ is the starting point. 
The space $\gh$ has to be replaced by the subspace $\ghb$ of
bounded holomorphic sections.
An orthonormal basis of the subspace defines a map
\begin{equation}
M\to \P(\ghb)\ .
\end{equation}
This defines an embedding into the infinite dimensional
projective space.
By the continuity of the evaluation functional
\refE{qhdef} Riesz's theorem can also be applied to define
the coherent vectors. For more details see \cite{Raw}, \cite{Wood}, \cite{Od},
 \cite{Odcs}.

\medskip

We need also the \emph{coherent projectors} used by Rawnsley
\begin{equation}\label{E:cohproj}
P_{\pi(q)}=\frac {|e_q\rangle\langle e_q|}{\langle e_q,e_q\rangle}\ .
\end{equation}
Here we used the convenient bra-ket notation.
For $s,t\in \gh$  the symbol $|s\rangle\langle t|$ 
denotes the following rank 1 operator of $\gh$ (resp. of  $\gul$)
\begin{equation}
|s\rangle\langle t|: r\to \langle t,r\rangle \cdot s
\end{equation}
By the normalization the projectors are 
indeed only depending on the points $\pi(q)$ of the manifold.

Rawnsley introduced the Epsilon function
\begin{equation}\label{E:epsilon}
\epsilon(\pi(q)):=|q|^2\langle e_q,e_q\rangle ,\quad
\text{with}\quad
|q|^2:=h(\pi(q))(q,q).
\end{equation}
Let $s_1$ and $s_2$ be two sections. At a fixed point 
$x=\pi(q)$ we can write $s_1(x)=\hat q(s_1)q$ and 
$s_2(x)=\hat q(s_2)q$ and hence using \refE{qhdefa} 
\begin{equation}
h(s_1,s_2)(x)=
\overline{\hat q(s_1)}\cdot \hat q(s_2)\cdot |q|^2
=
\langle s_1,e_q\rangle \langle e_q,s_2\rangle |q|^2=
\langle s_1,P_xs_2\rangle \cdot \epsilon(x)\ .
\end{equation}
After integration we obtain the over-completeness property
of the coherent states
\begin{equation}\label{E:cohover}
\langle s_1,s_2\rangle =\int_M
\langle s_1,P_xs_2\rangle\epsilon(x) \Omega(x)\ .
\end{equation}
We calculate
\begin{equation}\label{E:epssec}
\epsilon(x)=|q|^2\langle e_q,e_q\rangle
=|q|^2\sum_{j=0}^N|\hat q(s_j)|^2=
\sum_{j=0}^N|\hat q(s_j)|^2 h(x)(q,q)=
\sum_{j=0}^N h(s_j,s_j)(x)\ .
\end{equation}
It was shown in \cite[Equ. (3.4)]{CGR1} that
 for $\epsilon\equiv const$ one obtains
\begin{equation}\label{E:epsconst}
\epsilon=\frac {\dim\gh}{\volu (M)}\ .
\end{equation}
On homogeneous K\"ahler manifolds with a homogeneous quantum line bundle
(in particular also with homogeneous metric)
the function  $\epsilon(x)$ is invariant under moving the point,
hence it is a constant.
In particular, it is constant for  the projective space $\pnc$.
See \refP{projind} for more information.

\medskip
To compare this approach with the local description used by Berezin
we have to choose a section
$s_0\in\gh$, $s_0\not\equiv 0$. Let  
$V=\{x\in M\mid s_0(x)\ne 0\}$ be the open
subset
on which the section does not vanishes\footnote
{In the terminology of \refS{polar}
we remove the support of the divisor of the section $s_0$.
}.
Now $s_0$ is a holomorphic frame for the bundle $L$ over $V$. This says that 
over $V$ every holomorphic (differentiable) section
can be described as $s(x)=\hat s(x)s_0(x)$
with a holomorphic (resp. differentiable) function $\hat s$. The mapping 
$s\mapsto\hat s$ defines an isometry of  $\gh$ (resp. of $\gul$)
 into the  L${}^2$
space of holomorphic (resp. differentiable)
 functions on $V$ with respect to the measure
$\ \mu_{s_0}(x)=h(s_0,s_0)(x)\Omega(x)$.

With respect to the frame $s_0$ the function $\hat h$ describing the metric
is given as $\hat h(x)=h(s_0,s_0)(x)$. Hence we can describe 
the  scalar product for $\varphi,\psi\in\gul$ as
\begin{equation}\label{E:skpber}
\langle \varphi,\psi\rangle=
\int_V\overline{\hat \varphi(x)}\hat\psi(x)h(s_0,s_0)(x)
\Omega(x)\ .
\end{equation}
If  we introduced the local K\"ahler potential 
given by \refE{kaehlerpot}.
then this can be rewritten as
\begin{equation}\label{E:skpberb}
\langle \varphi,\psi\rangle=
\int_V\overline{\hat \varphi(x)}\hat\psi(x)\exp(-K(x))
\Omega(x)\ .
\end{equation}
It is enough to calculate the integral on $V$, because $M\setminus V$
is of (complex) codimension 1, hence of measure zero, see \refS{polar}.
Such a description is always possible. For doing explicit calculations
Berezin considered special cases where $V$ is either $\C^n$, or 
a subset of special type of $\C^n$ 
(e.g. bounded symmetric domains) 
\cite{Berebd},\cite{Berequ}, \cite{Beress}.

\section{The polar divisor}\label{S:polar}
\subsection{The definition of the polar divisor}\label{SS:polar}
\begin{definition}\label{D:polar}
Let $x\in M$ and $0\ne q\in\pi^{-1}(x)$ then the 
\emph{polar divisor} $\S_x$ associated to $x\in M$ is defined as
\begin{equation}
\S_x:=\{\, x'\in M\mid \langle e_q,e_{q'}\rangle=0\quad\text{for\  }
0\ne q'\in\pi^{-1}(x')\;\}\ .
\end{equation}
I.e. the polar divisor $\S_x$ is the set of points
on the manifold for which the associated coherent vectors
are orthogonal to the coherent vectors associated to $x$.
\end{definition}
Due to the relation \refE{cohtrans} the definition is independent of the
representing elements $q$ and $q'$.
In the context of Perelomov's coherent states the notion of 
 polar divisors  was introduced by Berceanu
 \cite{Berpre94},\cite{BerJGP},\cite{Berbia96},\cite{Berbia97}.
On the purely geometric side the polar divisor was used earlier by 
H.H.~Wu \cite{Wu} for the  complex Grassmannians.
As we will see in the following it has a meaning in much more general
situations.
There should be no danger of confusion with the notion of 
polar divisor in complex analysis as the divisor of the
polar part of a meromorphic function. (See the remark after
Equation \refE{come} for a connection.)

Note that for every (meromorphic or holomorphic) section of a line bundle
there is an associated divisor in the sense
of algebraic geometry. 
For a thorough treatment of the relation between divisors, line bundles and
sections of line bundles see \cite[p.130ff]{GH}, \cite{SchlRS}.
What we need here are only the following facts.
For a holomorphic section $s\not\equiv 0$ of a line bundle the zero-set of 
the section   can be
decomposed into a union of 
(complex) one-codimensional ``irreducible subvarieties'' which are
not necessarily smooth.
The complement of the zero-set is an open dense subset of
 $M$\footnote{We assume $M$ to be connected and compact}. 
 Each codimension one {\it irreducible subvariety}
can be given locally as zero-set of an algebraic function.
By the irreducibility the vanishing order along the subvariety is 
constant. Hence we can assign to the section $s$ the formal sum $(s)$
of (irreducible) codimension one subvarieties with integer coefficients  
\begin{equation}\label{E:divs}
(s):=\sum_{Y\text{ irreducible}\atop{\text{subvarity of } M
\atop \text{of codimension } 1}}
 n_Y Y\ ,
\end{equation}
where $n_Y$ denotes the vanishing order along $Y$.
By the compactness of $M$ the sum \refE{divs} will always be finite.
Every such  formal sum with $n_Y\in\Z$ fulfilling the restriction 
that $n_Y\ne 0$ only for finitely many $ Y $ 
is called a {\it divisor} of $M$.
The sum $(s)$ is called the {\it divisor of the section} $s$.
For meromorphic sections negative integers (corresponding to 
algebraic poles) are allowed.
Two divisors are called {\it linearly equivalent} if their
difference (as formal sum) is the divisor of a meromorphic
function on $M$.
Note that the functions are the (meromorphic) sections of the trivial 
line bundle.
By this an equivalence relation is defined.
The linear equivalence class of a divisor is called a divisor class.
The set of divisor classes carries a natural structure of
an abelian group under  addition of divisors.
For smooth projective varieties (as $M$ is one by the 
quantization condition) this
divisor class group is isomorphic to the group of 
isomorphy classes of algebraic line bundles, 
where for the latter the
group structure is defined by the tensor product of line bundles.
The isomorphism is given by assigning to the line bundle 
the divisor class of any non-trivial meromorphic section.
Note that the divisors of two meromorphic sections of the same  
line bundle are linearly equivalent.

Recall that the coherent vector  $e_q$ is a section of the
quantum line bundle and 
as introduced above  its 
divisor is given by $(e_q)$. 
The zero-set  of  $e_q$ 
(forgetting the multiplicities) is called the   
the \emph{reduced support} $red(e_q)$ of the divisor $(e_q)$.
\begin{proposition}\label{P:polar}
The polar divisor associated to $x$ is the reduced support of a divisor.
More precisely,
\begin{equation}
\S_x=red(e_q),\quad\text{with}\quad  q\in\pi^{-1}(x),\ q\ne 0\ .
\end{equation}
\end{proposition}
\begin{proof}
Let $q'\in L_0$, resp. $x'\in M$, $\pi(q')=x'$.
{}From \refE{qhdefa} it follows
$\ \langle e_{q'},e_q\rangle=\widehat{q'}(e_q)$
and 
\begin{equation}\label{E:eqq}
e_q(x')=e_{q}(\pi(q'))=\widehat{q'}(e_q)\cdot q'=
\langle e_{q'},e_q\rangle\cdot q'\ .
\end{equation}
Hence, $x'=\pi(q')$ is a zero of the section $e_q$ if and only if
$\langle e_{q'},e_q\rangle=0$. This shows the claim.
\end{proof}
By the above proof we  see that the multiplicity structure
of the zeros of $e_q$ and that of $\langle e_{q'},e_q\rangle=0$
are the same.  Hence we can indeed consider  $\S_x$ as  a divisor
if we assign to it the corresponding multiplicity of its components.
Note that due to relation \refE{cohtrans} 
\begin{equation}
(e_q)=(e_{q'})\quad\text{for}\quad
q,q'\in\pi^{-1}(x)\setminus\{0\},\quad x\in M\ .
\end{equation}
Hence we can assign for every $x\in M$ the divisor
\begin{equation}\label{E:divcs}
(e_x):=(e_{\pi^{-1}(x)})\ , 
\end{equation}
to the coherent state $e_x$.

We obtain
\begin{corollary}\label{C:zero}
The polar divisor $\S_x$ associated to $x$ is the divisor $(e_x)$
in the sense of algebraic geometry \refE{divcs}
  of the coherent state  $e_x$ associated to $x$.
\end{corollary}
Let $V$ be an open non-empty subset over which 
the bundle $L$ can be (holomorphically) trivialized, i.e.
$L_{|V}\cong V\times\C$. 
We take
$q'$ above $x'\in V$ as  $x'\mapsto (x',1)$
(i.e. we take as $q'$ the value at $x'$  of the frame given by the 
trivialization) 
then using  \refE{eqq} we obtain
\begin{proposition}\label{P:locrep}
The function 
$\ \hat e_q(x'):=\langle e_{(x',1)},e_q\rangle\ $
 is the holomorphic local representing
function for the  section  $e_q$.
\end{proposition}
Let $s$ be a global holomorphic section not identically zero.
The complement of the divisor $(s)$ 
\begin{equation}
  \label{E:Vs}
  V_s:=V\setminus (s)=\{x\in M\mid s(x)\ne 0\}
\end{equation}
will be an open dense subset of $M$.
If we apply \refP{locrep} to the trivialization obtained by taking as frame the
holomorphic section $s$
on $V_s$ then $(x',1)\cong s(x')$ and we can reformulate
\refP{locrep} as 
\begin{corollary}\label{C:locrep}
Let $s\not\equiv 0$ be a global holomorphic section. Then with
respect to the frame  given by the section $s$  
a local representing function over $V_s$ 
for the coherent vector $e_q$ is 
given by
\begin{equation}\label{E:locsec}
\hat e_q(x')=\langle e_{s(x')},e_q\rangle\ .
\end{equation}
\end{corollary}
Immediately from the definition of polar divisors we get
\begin{proposition}\label{P:polinvar}
The polar divisors obey the symmetry relation
\begin{equation}\label{E:polsymm}
y\in \S_x  \quad \longleftrightarrow\quad
x\in \S_y \ .
\end{equation}
\end{proposition}

\begin{remark}
By the definition \refE{qhdef} of $\hat q$, resp. of $e_q$
\refE{qhdefa}
one concludes that $\langle e_{q'},e_{q}\rangle$ varies
antiholomorphically in $q$ and
holomorphically in $q'$. Hence $q\mapsto \S_q$ defines an
antiholomorphic
family of divisors on $M$ and $q\mapsto\hat e_q$ an antiholomorphic
family of sections for the bundle $L$.
\end{remark}

\begin{remark}
By Bertini's theorem \cite{Hart}
the divisor of a generic global holomorphic section of the bundle $L$
is a smooth hypersurface. The divisors $(e_q)$ for the coherent
vectors are not necessarily generic, so one can not expect them
to be smooth in general.
See \cite{BerGr} for an example.
\end{remark}

\begin{example}
Let us consider the simplest example, the projective line $\pnc[1]$,
resp. the sphere $S^2$ with the K\"ahler structure given by 
$m$ ($m\in\N$) times the Fubini-Study form
\begin{equation}\label{E:p1}
\w=\frac {\i}{(1+z\zb)^2}dz\wedge d\zb
\end{equation}
as K\"ahler form with respect to the quasi-global coordinate $z$.
The corresponding quantum line bundle is $H^{\otimes m}$, where
$H$ is the hyperplane bundle.
The coherent vectors in the standard affine chart are given as
\begin{equation}\label{E:cohp}
e_{\psi(w)}(z)=\frac {m+1}{2\pi}(1+\wb z)^m\psi(z)\ ,
\end{equation}
where we take the (on this chart) non-vanishing section
$\psi$ as reference section.
The divisors are formal sums of points with integer coefficients.
We denote the divisor corresponding to the point with 
the coordinate $z_0$ by $\divi{z_0}$. 
The polar divisors are calculated  (using \refC{zero}) directly as the 
zero set of the section \refE{cohp}, hence
\begin{equation}
\S_w=m\divi{-\frac 1{\wb}},\quad w\ne 0,\infty,\qquad
\S_0=m\divi{\infty},\qquad \S_\infty=m\divi{0}\ .
\end{equation}
Using the original \refD{polar} we can also calculate
(compare \refE{locsec})
\begin{equation}\label{E:prod}
\langle e_{\psi (w')},e_{\psi (w)}\rangle_{H^{\otimes m}}
=
\frac{m+1}{2\pi}(1+{\overline{w}}w')^m\ ,
\end{equation}
yielding (in accordance with \refC{zero})
 clearly the same set of points where 
\refE{prod} vanishes.

In particular for $m>1$ the divisors appearing as divisors of coherent
states are not smooth because they have higher multiplicities.
\end{example}

The polar divisors appear at many places.
It was shown by Berceanu 
that for  Grassmannians and more general for Hermitian symmetric
spaces 
\cite{BerJGP},\cite{BerGr},\cite{Berbia96}
the polar divisor $\Sigma_x$ (with respect to Perelomov's
coherent states)  coincides with the cut-locus of the point $x$.
Recall that for a geodesic starting at $x$ the cut-point $y$ is
the point where the geodesic ceases to
be the shortest curve connecting $x$ and $y'$ with $y'>y$
on the geodesic. 
The cut-locus consists of all cut-points.
For more details see the above references.

As we will see in the following section 
for arbitrary compact K\"ahler manifolds 
the polar divisors 
describe
the zero-sets of two-point functions and more
general of $m$-point functions.
They  appear as singularity sets of the analytic extensions
of real-analytic metrics in the bundle (see \refE{comec})
and as  singularity sets of the covariant two-point
Berezin symbols.

\section{Cauchy formulas and multi-point functions}\label{S:cauchy}
\subsection{Coherent projective K\"ahler embedding}
{~}

\vspace*{0.2cm}

Let the compact K\"ahler manifold $(M,\w_M)$ be embedded via the
(holomorphic) coherent state map \refE{cohmapc} with respect to the 
very ample line bundle $(L,h_L)$:
\quad $i:M\hookrightarrow \pnc$. Fix the orthonormal sections $s_j$,
$j=0,\ldots,N$  
of the quantum bundle $L$.
Let $t_j,\ j=0,\ldots,N$ be the sections of the hyperplane bundle $H$ over
$\pnc$ corresponding to the linear forms $Z_j,\ j=0,\ldots,N$.
By construction we have $i^*H\cong L$ and $i^*(t_j)=s_j$
(i.e. $s_j(x)=t_j(i(x)$)
under the identification of the bundles.
It is well-known that the $t_j$ are orthogonal sections of $H$ with
norm independent of $j$.
We will denote the rescaled orthonormal section by $t_j'$ and obtain
$s_j=\tau i^*(t_j')$ with a factor 
 $\tau$ independent of $j$.
More precisely, $\tau=\sqrt{\volu(\pnc)/(N+1)}$.
Note that the pullback in our case is nothing else as the restriction
of the section to the embedded manifold $M$.

Consider the case where the coherent state embedding is an isometric
(projective)
K\"ahler embedding, i.e. the pullback of the Fubini-Study form $\w_{FS}$
coincides with $\w_M$. By \refP{invarmet}
the metric $h_L$ in the bundle $L$ is up to a positive 
scalar multiple the pullback of the
metric $h_{FS}$ in the hyperplane  bundle: $h_L=\rho \cdot i^*h_{FS}$, 
$\rho\in\R,\ \rho> 0$.
Let $\epsilon_M$ be the Epsilon function \refE{epsilon} for the manifold $M$
and the bundle $L$.
Due to the explicit  description \refE{epssec} 
of the Epsilon function $\epsilon_M$ 
it is up to a constant the restriction of $\epsilon_{\pnc}$
to the embedded points.
 The latter is constant, hence also
$\epsilon_M$. In more detail:
\begin{multline}
\epsilon_M(x)=\sum _{j=0}^N h_L(s_j,s_j)(x)=
\sum _{j=0}^N h_L(\tau i^*(t_j'),\tau i^*(t_j'))(x)=
\\
=\rho\tau^2
\sum _{j=0}^N h_H(t_j',t_j')(i(x))=
\rho\tau^2\epsilon_{\pnc}(i(x))=\rho\tau^2\tfrac {N+1}{\volu(\pnc)}=\rho\ .
\end{multline}
{}From this $\rho$ calculates to $\frac {N+1}{\volu(M)}$.

Rawnsley calls a quantization where all the data can be obtained
by  pulling back  the objects: bundle, forms, etc.
via the (holomorphic) coherent state map from the projective space
to the manifold $M$ \emph{projectively induced}.
Hence projectively induced quantizations have constant Epsilon
functions.
In fact the converse is also true:
\begin{proposition}\label{P:projind}
(Cahen, Gutt, Rawnsley, \cite[p.58]{CGR1})
A quantization of $(M,\w)$ with quantum line bundle $(L,h)$ is
projectively induced
if and only if the Epsilon function is constant.
\end{proposition}

This has  very interesting consequences for compact homogeneous 
K\"ahler manifolds $M\cong G/H$ with a homogeneous quantum line bundle.
Recall that it is assumed in this case 
that the K\"ahler form $\w_M$ and the metric
 in the bundle are invariant under
the action  of the group $G$.
In particular, $\epsilon_M$ will be constant, hence

\begin{corollary}
For  a compact homogeneous K\"ahler manifold which admits a homogeneous very
ample quantum line bundle $L$ the coherent state embedding using this
bundle is a
projective K\"ahler embedding.
In such cases the K\"ahler form is the pullback of the
Fubini-Study form.
\end{corollary}

\bigskip
\subsection{Cauchy formulas}\label{SS:cf}
{~}

\vspace*{0.2cm}

Let us return to the general situation of the coherent state 
embedding $i$, \refE{cohmapc} without assuming it to be a K\"ahlerian 
embedding.
Fix an orthonormal basis
 $s_0,s_1,\ldots,s_N$ of the sections of the quantum line
 bundle $L$. 
Denote by $\phi$ the map from $L_0$ to $\C^{N+1}$ defined by 
the composition 
 \begin{equation}\label{E:cohvmap}
q\quad \mapsto\quad  \overline{e_q}\quad \mapsto\quad
\phi(q):=({\hat q(s_0)},{\hat q(s_1)},...,
{\hat q(s_N)})\ .
\end{equation}
Clearly, $i(\pi(q))=[\phi(q)]$.

In the following,  three scalar products will  appear:
(1) $\langle .,.\rangle_{L}$, the scalar product on the space of
global sections of $L$ given by \refE{skp},
(2) $\langle .,.\rangle_{\C^{N+1}}$, the standard scalar product on
$\C^{N+1}$,
and (3)  $\langle .,.\rangle_{H}$
the scalar product on the space of global sections 
of the hyperplane bundle $H$ on $\pnc$.
Again (3) is defined by \refE{skp}, but now the manifold is $\pnc$ and 
the sections are the hyperplane sections.
Recall that all our scalar products are conjugate linear
 in the first arguments.
We will call relations between these scalar products
(evaluated for coherent vectors) Cauchy formulas.
The first Cauchy formula is immediate from 
\refE{cohdecomp}, \refE{cohvmap}
\begin{proposition}\label{P:cauchya}
In the above situation we have 
\begin{equation}\label{E:cauchya}
{\skp {e_{q}} {e_{q'}}}_L
={\skp {\phi({q'})} {\phi({q})}}_{\C^{N+1}}\ .
\end{equation}
\end{proposition}
Next we want to find  relations between the scalar product of 
coherent vectors 
of $L$ over $M$ 
and the scalar product of coherent vectors of the hyperplane bundle $H$
over the projective space.
Let $s\not\equiv0$ be a holomorphic section of $L$, which is non-vanishing
over the dense open  subset $V_s$.
Clearly, the scalar product of two coherent states is not defined.
But, after choosing such a section we can set
\begin{equation}\label{E:scs}
\skps {e_x}{e_y}{s}:=\skp {e_{s(x)}}{e_{s(y)}}\ .
\end{equation}
If we choose another holomorphic section $s'\not\equiv 0$ then on 
$V_s\cap V_{s'}$ 
we have $s'(x)=f(x)s(x)$ with $f$ a 
non-vanishing holomorphic function on
this set.
Hence,
\begin{equation}\label{E:scst}
\skp{e_x}{e_y}_{s'}=
\frac 1{f(x)}\frac 1{\overline{f(y)}}\;
\skp{e_x}{e_y}_{s}\ .
\end{equation}

Recall that Rawnsley's Epsilon function 
 $\epsilon$ \refE{epsilon} can be written as 
\begin{equation}
\epsilon_M(x)={|s(x)|}^2\skp{e_{s(x)}}{e_{s(x)}}\ ,
\quad\text{where}\quad
{|s(x)|}^2:=h(s(x),s(x))\ .
\end{equation}
The function 
$$
\chi_{s}(x,x):=\skp{e_{s(x)}}{e_{s(x)}}
={h(s(x),s(x))}^{-1}\cdot \epsilon_M(x)
$$
is real analytic and admits a real analytic extension to
the function 
\begin{equation}\label{E:come}
\chi_{s}(x,y):=\skp{e_x}{e_y}_{s}=
\skp{e_{s(x)}}{e_{s(y)}}\ ,
\end{equation}
which is holomorphic in $x$ and antiholomorphic in $y$.

\begin{remark}
Assume the metric $h$ in the bundle to be real-analytic.
Then from \refE{come} one concludes that for fixed $x$  the singularity 
set (in the variable $y$) of the extension of the metric 
is given by the polar divisor $(e_x)$.
\end{remark}
\begin{remark}
For $\epsilon(x)=\epsilon_M$ a constant, we see that
the scalar product of the coherent vectors is essentially given 
by the inverse of the local metric:
\begin{equation}\label{E:comec}
\skps {e_x}{e_y}{s}=\frac {\epsilon_M}{h(s(y),s(x))}\ .
\end{equation}
\end{remark}
\noindent
Let us apply this to $M=\pnc$ with the hyperplane bundle $H$ as 
quantum line bundle and the metric of the bundle induced by the
Fubini-Study metric.
Let 
\newline
$\ V_0:=\{[z]=(z_0:z_1:\cdots:z_N)\mid z_0\ne 0\}\ $
be the standard affine chart. The points of $V_0$ 
can be given in a normalized way
as $(1:w)$ with $w\in\C^N$.
Take $t_0$ to be the section of the 
hyperplane bundle corresponding to the linear form $Z_0$, i.e.
$t_0(w)=1$ for all $w$. We set
$\hat y:=(t_0(y),y)=(1,y)$ and $\hat x=(t_0(x),x)=(1,x)$.
Then
\begin{equation}
h(t_0(y),t_0(x))=\frac 1{1+\bar y\cdot x}=
\frac 1{{\skp {\hat y}{\hat x}}_{\C^{N+1}}}\ .
\end{equation}
Hence in this case \refE{comec} specializes to
\begin{equation}\label{E:exypn}
{\skp{e_x}{e_y}}_{t_0}=\skp{e_{t_0(x)}}{e_{t_0(y)}}
={\skp {\hat y}{\hat x}}_{\C^{N+1}}\cdot F\ ,
\end{equation}
with
\begin{equation}\label{E:Fdef}
F:=\epsilon_{\pnc}=\frac {N+1}{\volu(\pnc)}\ .
\end{equation}
Now we return to the general situation.
The second Cauchy formula is expressed in
\begin{proposition}\label{P:cauchyb}
Let $q,q'\in L_0$ with $\pi(q)=x$ and  $\pi(q')=y$, $i:M\to\pnc$ the
coherent state embedding \refE{cohmapc}, and $\phi$ the map \refE{cohvmap},
then on the affine chart $V_0$
\begin{equation}\label{E:cauchyb}
{\skp{e_q}{e_{q'}}}_L=
\widehat{q}(s_0)\overline{\widehat{q'}(s_0)}\cdot{\skp{\widehat{i(y)}}
{\widehat{i(x)}}}_{\C^{N+1}}=
\frac {\widehat{q}(s_0)\overline{\widehat{q'}(s_0)}}{F}\cdot
{\skp {e_{i(x)}}{e_{i(y)}}}_{t_0} \ .
\end{equation}
\end{proposition}
\begin{proof}
We start from \refE{cauchya} in \refP{cauchya} and divide
 the vectors $\phi(q)$ and $\phi(q')$ on the left
hand side by their first components.
 This can be done because we are on $V_0$.
We obtain
\begin{equation*}
{\skp{e_q}{e_{q'}}}_L= \widehat{q}(s_0)\overline{\widehat{q'}(s_0)}\;
{\skp{\phi(q')_{norm}}{\phi(q)_{norm}}}_{\C^{N+1}}\ .
\end{equation*}
Here $\phi(q)_{norm}$ is the normalized representative
which  has first component 1.
Using $\phi(q)_{norm}=\widehat{i(\pi(q))}$ and \refE{exypn} we obtain
\begin{equation*}
{\skp{e_q}{e_{q'}}}_L=\widehat{q}(s_0)\overline{\widehat{q'}(s_0)}\;
{\skp{\widehat{i(y)}}{\widehat{i(x)}}}_{\C^{N+1}}=
\frac {\widehat{q}(s_0)\overline{\widehat{q'}(s_0)}}{F}\;
{\skp {e_{i(x)}}{e_{i(y)}}}_{t_0}\ .
\end{equation*}
\end{proof}
The third Cauchy formula is expressed in 
\begin{theorem}\label{T:cauchy}
Let $(M,\omega)$ be a quantizable K\"ahler manifold with 
very ample quantum line bundle
$L$. Let $i:M\to\pnc$ be the  coherent state
embedding \refE{cohmapc}, $H$ the hyperplane section bundle.
For every section $t$ of $H$ denote 
by $i^*(t)$ its pullback to $M$. Assume $t\not\equiv 0$  then
over $V_t:=\{z\in\pnc\mid t(z)\ne 0\}$
\begin{equation}\label{E:cauchy}
\skp {e_{x}}{e_{y}}_{i^*(t)}=
\frac {\volu(\pnc)}{N+1}\;{\skp {e_{i(x)}}{e_{i(y)}}}_{t}\ .
\end{equation}
\end{theorem}
\begin{proof}
First consider the section $t_0$. In this case $i^*(t_0)=s_0$.
Note that in view of \refE{cauchyb} and \refE{Fdef} it is
enough to show 
that $\ \widehat{s_0(x)}(s_0)=1\ $.
But by definition 
$\ \widehat{s_0(x)}(s_0)\cdot s_0(x)=s_0(x)$, hence 
\begin{equation}\label{E:cauchys}
{\skp {e_{x}}{e_{y}}}_{s_0}=
{\skp {e_{s_0(x)}}{e_{s_0(y)}}}=
\tfrac 1F {\skp {e_{i(x)}}{e_{i(y)}}}_{t_0}\ .
\end{equation}
\newline
Now take a  general $t\not\equiv 0$. Recall that the complement 
of a zero-set of a section ($\not\equiv 0$) is always a dense open subset.
Hence the same is true for finite intersections of such sets.
On the dense open set $V_0\cap V_t$
we have $t(z)=f(z)\cdot t_0(z)$
with  $f(z)$ a holomorphic function on the intersection.
For the pull-backs we obtain
\begin{equation*}
(i^*t)(x)=t(i(x))=f(i(x))\cdot t_0(i(x))=(i^*f)(x)\cdot (i^*t_0)(x)\ .
\end{equation*}
This implies (using \refE{scst})
\begin{gather*}
\skp {e_{x}}{e_{y}}_{i^*(t)}
=((i^*f)(x))^{-1}\cdot (\overline{(i^*f)(y)})^{-1}\cdot
\skp {e_{x}}{e_{y}}_{i^*(t_0)}, 
\\
{\skp {e_{i(x)}}{e_{i(y)}}}_{f\cdot t_0}
=(f(i(x)))^{-1}\cdot(\overline{f(i(y))})^{-1}
{\skp {e_{i(x)}}{e_{i(y)}}}_{t_0}\ .
\end{gather*}
But note that $i^*(t_0)=s_0$ and $f(i(x))=i^*f(x)$. 
The claim
follows from \refE{cauchys}.
\end{proof}

The reason for calling the 
Equations  \refE{cauchya}, \refE{cauchyb}, \refE{cauchy}
``Cauchy formulas'' is that in the case
of the Grassmannians the appearing formulas are 
essentially the (Binet-) Cauchy formulas \cite[p.10]{Gant}
which give relations between the intrinsic metric on the
Grassmannian and the pull-back of the Fubini-Study metric 
obtained via the Pl\"ucker embedding, see
 \cite[Equations 3.13, 4.7]{BerJGP},\cite{BerGr}.

\bigskip

\subsection{Two-point functions}
{~}

\vspace*{0.2cm}
As already noted in \refSS{polar} the assignment
\begin{equation}\label{E:2pvf}
L_0\times L_0\quad\to \quad\C,\qquad
(q,q')\mapsto \skp {e_{q}}{e_{q'}}
\end{equation}
defines a real-analytic function holomorphic in $q$ and
 antiholomorphic in $q'$. It can be normalized by setting
\begin{equation}\label{E:2pf}
\phi(q,q'):=\frac {\skp {e_{q}}{e_{q'}}}{\Vert e_{q}\Vert\Vert
e_{q'}\Vert }\ .
\end{equation}
By the Cauchy-Schwartz inequality its absolute value is bounded by 1.

Using \refE{cohtrans} we see
\begin{equation}
\phi(cq,c'q')=\frac {|c|}{c} \frac {|c'|}{\overline{c'}}\,
\phi(q,q'),\qquad c,c'\in\C^*\ .
\end{equation}
Due to the appearing phase factors,
it does not descend to a two-point function on $M$.
Clearly, one way out is to take the modulus  (or its square,
see \cite{CGR2})
of \refE{2pf}. In this way one obtains the function 
$\ \psi:M\times M\to[0,1]$ 
\begin{equation}\label{E:2pfm}
\psi(x,y):=\frac {{|\skp {e_{q}}{e_{q'}}|}^2}{{\Vert e_{q}\Vert}^2{\Vert
e_{q'}\Vert}^2 }, \qquad x=\pi(q),\ y=\pi(q')\ .
\end{equation}
This is a globally defined real-valued real-analytic function.
Unfortunately, the information contained in the complex phase gets lost.

Note that the set of zeros $M(\psi)$ of \refE{2pfm} can be given with the help
of polar divisors
\begin{equation}
M({\psi}):=\{\ (x,y)\in M\times M\mid y\in\S_x\ \}\ .
\end{equation}
Recall that by \refP{polinvar} the condition is symmetric in $x$ and $y$.
Clearly, the zero set $M(\phi)$ of
 \refE{2pf} consists of the fibers over $M(\psi)$.

\medskip
Take a section $s\not\equiv 0$ of $\gh$ and choose it as frame over
$V_s=M\setminus (s)$. We define
the function 
\begin{equation}\label{E:2pfs}
\hat\phi_s(x,y):=\frac {{\skp {e_{s(x)}}{e_{s(y)}}}}
{{\Vert e_{s(x)}\Vert}{\Vert
e_{s(y)}\Vert} }
\end{equation}
on  $V_s\times V_s$.
It ``represents'' the two-point function \refE{2pf}. But
note that  $\hat\phi_s$ depends on the section $s$ 
which was chosen as frame.

Immediately from the definition we get
\begin{equation}
 \hat\phi_s(x,y)=\overline{\hat\phi_s(y,x)}\ . 
\end{equation}
If we deal with different manifolds and if there is a danger
of confusion, we will exhibit also the manifold in the notation
of the two-point function (resp. of the  $m$-point functions
introduced later on).

\begin{proposition}
\label{P:2pf}
Let $i$ be the coherent state embedding \refE{cohmapc}.
Let $t\not\equiv 0$ be a section of the hyperplane section bundle $H$ 
and $s=i^*(t)$ the corresponding section of the quantum line bundle
then

\noindent
(a) 
\begin{equation}\label{E:2pfis}
\hat\phi_{M,\,i^*(t)}(x,y)=
\hat\phi_{\pnc,\,t}(i(x),i(y))
=i^*(\hat\phi_{\pnc,\,t})(x,y)\ .
\end{equation}
(b) For the first coordinate function $s_0=i^*(t_0)$ and with
$\widehat{i(x)}\in\C^{N+1}$ the normalized
homogeneous representative of $i(x)$ one has
\begin{equation}\label{E:2pfn}
\hat\phi_{M,s_0}(x,y)=
\frac {\skp {\widehat{i(y)}}{\widehat{i(x)}}}
{\Vert\widehat{i(x)}\Vert \Vert\widehat{i(y)}\Vert}\ .
\end{equation}
(c)
\begin{equation}\label{E:2fm}
\psi_{M}(x,y)=i^*\psi_{\pnc}(x,y)=
\psi_{\pnc}(i(x),i(y))
=
\frac {{|\skp {\widehat{i(y)}}{\widehat{i(x)}}|}^2}
{||\widehat{i(x)}||^2 ||\widehat{i(y)}||^2}\ .
\end{equation}
\end{proposition}
\begin{proof}
(a) is immediate from   
\refT{cauchy}.
\newline
(b) follows from (a) using \refE{exypn}.
\newline
(c) follows from (a), resp. (b) by taking the
squared modulus. Note that
$\psi_M=\hat\phi_{M,\,s}\cdot \overline{\hat\phi_{M,\,s}}$
independently on the section $s$ chosen.
\end{proof}
In the case of Perelomov's coherent 
states the Equation \refE{2pfn} was also 
called a Cauchy formula in \cite[Equ. 3.13]{BerJGP}. 

The two-point functions (complex-valued or real-valued)
play an important r\^ole.
{}From their very definition they give the transition amplitudes
for coherent states.
They appear as  integral kernel of the
Berezin transform 
which relates contravariant and covariant Berezin
symbols, see \cite{Schlbia98}.
See also the discussion in \cite{CGR2} for the real-valued
two-point function and its relation to Calabi's diastatic function $D$.
Let us add a few words on this relation.
For real-analytic 
metrics $h$ also another two-point function
$\tilde\psi$ is introduced in the article \cite{CGR2} . 
It is given completely in local terms of the metric.
The relation $\tilde\psi=\exp(-D/2)$ is shown.
Certain natural behaviour under pull-backs is proven.
In the case that $\epsilon=const$ (in the terminology of 
\cite{CGR2}: the bundle is regular) one obtains $\tilde\psi=\psi$.
The key ingredients for this is equation \refE{comec} 
which relates the global scalar product with the local metric.
For regular line bundles \refE{2fm} was also 
proven in \cite{CGR2}.


\bigskip
\subsection{Cyclic $m$-point functions and the three-point function}
\label{SS:3pf}
{~}

\vspace*{0.2cm}
Let us consider the (cyclic) $m$-point function for $m\in\N, \ m\ge 2$
\begin{equation}\label{E:mpf}
\begin{gathered}
\Psi^{(m)}:
M\times M\times.....\times M\quad\to \quad\{z\in\C\mid |z|\le 1\},\\
\Psi^{(m)}(x^{(1)},x^{(2)},\ldots,x^{(m)})=
\frac {\skp {e_{q^{(1)}}}{e_{q^{(2)}}}\skp {e_{q^{(2)}}}{e_{q^{(3)}}}
\ldots 
\skp {e_{q^{(m)}}}{e_{q^{(1)}}}}
{{\Vert e_{q^{(1)}}\Vert}^2{\Vert e_{q^{(2)}}\Vert}^2\cdots
{\Vert e_{q^{(m)}}\Vert}^2},\\
 x^{(i)}=\pi(q^{(i)}),\ i=1,\ldots ,m\ .  
\end{gathered}
\end{equation}
It is a complex-valued and real-analytic function in its variables.
Note that the phase ambiguity of the lifts is  canceled by
this combination.
The function $\Psi^{(m)}$ 
can be written in terms of the complex-valued two-point function
as
\begin{equation}\label{E:mpfp}
\Psi^{(m)}(x^{(1)},x^{(2)},\ldots,x^{(m)})=
\hat\phi_s(x^{(1)},x^{(2)})\cdot
\hat\phi_s(x^{(2)},x^{(3)})\cdots
\hat\phi_s(x^{(m)},x^{(1)})\ 
\end{equation}
with respect to any section $s\not\equiv 0$ of $L$.
Note that  $\Psi^{(2)}=\psi$, the real-valued two-point function
as defined  in \refE{2pfm}.
But for $m>2$ the $\Psi^{(m)}$ will be complex-valued.
\begin{proposition}\label{P:mpfr}
Let $\Psi_M^{(m)}$, resp. $\Psi_{\pnc}^{(m)}$ be the $m$-point function of
the manifold $M$, resp. of the projective space. Let $i$ be the 
coherent state embedding \refE{cohmapc} and $\widehat{i(x)}$
an arbitrary homogeneous representative for the point $i(x)$ 
then 
\begin{equation}\label{E:mpfr}
\begin{gathered}
\Psi_M^{(m)}(x^{(1)},x^{(2)},\ldots,x^{(m)})=
i^*(\Psi_{\pnc}^{(m)})(x^{(1)},x^{(2)},\ldots,x^{(m)})=
\\
\Psi_{\pnc}^{(m)}(i(x^{(1)}),i(x^{(2)}),\ldots,i(x^{(m)}))
=\frac {\skp {\widehat{i(x^{(2)})}}{\widehat{i(x^{(1)})}}
\skp {\widehat{i(x^{(3)})}}{\widehat{i(x^{(2)})}}\cdots
\skp {\widehat{i(x^{(1)})}}{\widehat{i(x^{(m)})}}}
{||\widehat{i(x^{(1)})}||^2 ||\widehat{i(x^{(2)})}||^2\cdots 
 ||\widehat{i(x^{(m)})}||^2}\ .
\end{gathered}
\end{equation}
The function $\Psi_M^{(m)}$ is invariant under  cyclic 
permutations of its arguments. 
\end{proposition}
\begin{proof}
Using \refE{mpfp} we see that from \refE{2pfis} the first equality
follows. Now using \refE{2pfn} we obtain the last
equality.
Note that as $\widehat{i(x)}$ any homogeneous representative
can be chosen (but then it has be kept fixed).
The ambiguity will  cancel in this combination.
The invariance under cyclic permutations is clear.
\end{proof}
We can represent the last expression  in \refE{mpfr} as 
complex-conjugate of the similar expression where the $x^{(i)}$ 
appear in  strictly increasing index order $\mod m$.

Again, the zero-set of the $m$-point function can be given
with the help of  polar divisors.

Note that the last expression in \refE{mpfr} can be 
rewritten as follows.
Let $[u],[v]\in\pnc$ be points with homogeneous coordinates
$u,v\in\C^{N+1}$. The Cayley distance $d_C$  of the
two points, i.e the geodesic
distance with respect to the Fubini-Study metric  is given as
\begin{equation}
  \label{E:Cay}
  d_C([u],[v])=\arccos\frac {|\skp {u}{v}|}{\Vert u\Vert\Vert v\Vert}\ .
\end{equation}
Hence the last expression in \refE{mpfr} can be given 
in terms of the Cayley distances 
of neighbouring points and  an additional global phase factor 
$\Phi$ depending on the points.

Let us study the three-point function  $\Psi=\Psi^{(3)}$
in more detail. Take $[u],[v],[w]\in\pnc$ and let
$a=d_C([u],[v]), b=d_C([v],[w]), c=d_C([w],[u])$
be the Cayley distances. 
To avoid degenerate situations
assume that $0<a,b,c< \pi/2$.
In particular, no point should lie in the support of the polar divisors of
the other two.
We can write 
\begin{equation}
  \label{E:geo}
\frac {\skp uv\skp vw \skp wu}
{||u||^2 ||v||^2 ||w||^2}=\cos a\cdot \cos b\cdot\cos c \cdot \e^{-\i\Phi}\ ,
\end{equation}
with the phase factor $\Phi=\Phi(u,v,w)$ defined by this formula.
The phase factor is related to the shape invariant $\rho$ 
introduced by Blaschke and Terheggen \cite{BlTe} for $\pnc[2]$, 
resp. by Brehm \cite{Brehm} for $\pnc$:
 $\rho= \cos a \cos b\cos c \cos\Phi$.

Hangan and Masala showed that $\Phi$ has the following geometric meaning
\cite{HaMa}:
Take the (oriented) geodesic triangle $\sigma([u],[v],[w])$ with the
vertices  $[u],[v],[w]$, i.e. the 
surface 
swept out by the the geodesics between  the point $[w]$
and all points on the geodesic between $[u]$ and $[v]$.
Then
\begin{equation}
  \label{E:HaMa}
  \Phi=\int_{\sigma([u],[v],[w])}\w_{\pnc} +2k\pi,\quad k\in\Z.
\end{equation}
Recall that  $\w_{\pnc}$ is the Fubini-Study K\"ahler form of $\pnc$.
See also \cite{Berbia98} for a different proof by coherent
state methods.
Clearly, due to the fact that $\w_{\pnc}$ is closed the integral does not
change if we replace $\sigma([u],[v],[w])$ by any deformed surface
as long as the boundary is fixed.

By applying  \refP{mpfr} we see 
 that \refE{geo} is the complex conjugate of the three-point function of
the projective space. We obtain
\begin{theorem}\label{T:geo}
Let $(M,\omega)$ be a quantizable K\"ahler manifold with 
very ample quantum line bundle
$L$. Let $i:M\to\pnc$ be the  coherent state
embedding \refE{cohmapc}.
Then the three-point function 
\begin{equation}
\label{3pf}
\Psi^{(3)}(x,y,z)=
\frac {\skp {e_{q}}{e_{q'}}\skp {e_{q'}}{e_{q''}}
\skp {e_{q''}}{e_{q}}}
{{\Vert e_{q}\Vert}^2{\Vert e_{q'}\Vert}^2{\Vert e_{q''}\Vert}^2},\qquad
 x=\pi(q),\ y=\pi(q'),\ z=\pi(q'')
\end{equation}
can be written as 
\begin{equation}
\Psi_M^{(3)}(x,y,z)=\cos a\cdot\cos b\cdot\cos c\cdot\e^{\i\Phi}
\end{equation}
with 
$a=d_C(i(x),i(y))$, $b=d_C(i(y),i(z))$, $c=d_C(i(z),i(x))$
the Cayley distances in $\pnc$ and phase
\begin{equation}
\Phi=\int_{\tilde\sigma(i(x),i(y),i(z))}\w_{\pnc}
\end{equation}
where $\tilde\sigma(i(x),i(y),i(z))$ is any deformation of the geodesic
triangle (in $\pnc$) with fixed boundary given by the
geodesics (in  $\pnc$) connecting  the points $i(x),i(y)$, and $i(z)$.
\end{theorem}
The relationship between the phase $\Phi$ of the 3-point function
 and the symplectic area of a geodesic triangle on
the manifold itself
is studied for the complex Grassmann manifolds in 
\cite{Berbia98}.
\bigskip
\subsection*{Acknowledgements}
The authors gratefully acknowledge the Volkswagen-Stiftung
for their support in the Research-in-Pairs program.
We would also like
to thank the Mathematisches Forschungsinstitut in Oberwolfach 
and its staff for
the experienced  hospitality.
The work was finished during an invitation  of  S.B. at the
University of Mannheim.
S.B. likes to thank the DFG for supporting his stay.

\bibliographystyle{amsplain}
\ifx\undefined\bysame
\newcommand{\bysame}{\leavevmode\hbox to3em{\hrulefill}\,}
\fi


\end{document}